\documentclass{gtart_h}  

\def\ifplaintex{\expandafter\ifx\csname documentclass\endcsname\relax}

\ifplaintex 
\else
\fi

\expandafter\ifx\csname beginpicture\endcsname\relax
\expandafter\ifx\csname documentclass\endcsname\relax
\input pictex \else
\input prepictex \input pictex \input postpictex \fi\fi

\def\gt{{\mathsurround=0pt\it $\cal G\mskip-2mu$eometry \&\ 
$\cal T\!\!$opology}}        

\def\gtp{{\mathsurround=0pt\it $\cal G\mskip-2mu$eometry \&\ 
$\cal T\!\!$opology $\cal P\!$ublications}}  

\def\lognumber#1{\def\thelognumber{#1}}
\def\volumenumber#1{\def\thevolumenumber{#1}}
\def\papernumber#1{\def\thepapernumber{#1}}
\def\volumeyear#1{\def\thevolumeyear{#1}}

\def\pagenumbers#1#2{\def\startpage{#1}\def\finishpage{#2}}
\def\published#1{\def\publishdate{#1}}
\def\proposed#1{\def\theproposer{#1}}
\def\seconded#1{\def\theseconders{#1}}
\def\received#1{\def\receiveddate{#1}}
\def\revised#1{\def\reviseddate{#1}}
\def\accepted#1{\def\accepteddate{#1}}

\def\asciiaddress#1{\def\theasciiaddress{#1}}
\def\asciiemail#1{\def\theasciiemail{#1}}

\long\def\asciiabstract#1{\long\def\theasciiabstract{#1}}

\let\\\par\let\thelognumber\relax
\let\thevolumenumber\relax\let\thepapernumber\relax
\let\thevolumeyear\relax\let\thesamplenumber\relax\let\startpage\relax
\let\finishpage\relax\let\publishdate\relax\let\receiveddate\relax
\let\reviseddate\relax\let\accepteddate\relax\let\theasciititle\relax
\let\theasciiauthors\relax\let\theasciiaddress\relax
\let\theasciiabstract\relax
\let\theasciiemail\relax\let\theshortauthors\relax\let\theshorttitle\relax

\long\def\maketitlep{   

\count0=\startpage

\gt\hfill      
\beginpicture
\setcoordinatesystem units <0.33truein, 0.33truein> point at 2.2 0.9
\setplotsymbol ({$\cal G$})
\plotsymbolspacing=9truept
\circulararc 315 degrees from 0 1 center at 0 0
\setplotsymbol ({$\cal T$})
\circulararc 315 degrees from 1 -1 center at 1 0
\endpicture
\break
{\small\ifx\thesamplenumber\relax 
Volume \else Sample
\fi\thevolumenumber\ (\thevolumeyear)
\startpage--\finishpage\nl
Published: \publishdate}
\vglue 0.5truein plus 0.4fil minus 0.1truein

{\parskip=0pt\leftskip 0pt plus 1fil\def\\{\par\smallskip}{\ifplaintex\large
\else\Large\fi\bf\thetitle}\par\medskip}   

\vglue 0pt plus 0.1fil 

{\parskip=0pt\leftskip 0pt plus 1fil\def\\{\par}{\sc\theauthors}
\par\medskip}

\vglue 0pt plus 0.1fil 

{\small\parskip=0pt\let\newline\\
{\leftskip 0pt plus 1fil\def\\{\par}{\sl\theaddress}\par}
\expandafter\ifx\theemail\relax    
\relax\else\vglue 5pt plus 0.02fil minus 2pt\def\\{\stdspace{\rm 
and}\stdspace} 
\cl{Email:\stdspace\tt\theemail}\fi
\ifx\theurl\relax                  
\relax\else\vglue 5pt plus 0.02fil minus 2pt\def\\{\stdspace{\rm 
and}\stdspace}
\cl{URL:\stdspace\tt\theurl}\fi\par}

\vglue 7pt plus 0.3fil minus 3pt

{\bf Abstract}
\vglue 5pt plus 0.1fil minus 2pt

\theabstract

\vglue 7pt plus 0.3fil minus 3pt

{\bf AMS Classification numbers}\quad Primary:\quad \theprimaryclass

Secondary:\quad \thesecondaryclass

\vglue 5pt plus 0.3fil minus 2pt

{\bf Keywords:}\quad \thekeywords

\vglue 10pt plus 0.5fil minus 5pt

{\small  Proposed: \theproposer\hfill Received: \receiveddate\nl
Seconded: \theseconders\hfill 
\ifx\reviseddate\relax                         
Accepted: \accepteddate                        
\else
Revised: \reviseddate                          
\fi}
\eject
}       

\let\maketitlepage\maketitlep
\let\maketitle\maketitlepage

\font\phead=cmsl9 scaled 950
\font=cmsl9 scaled 1050
\font\pnum=cmbx10 scaled 913
\font\lnum=cmbx10 
\font\pfoot=cmsl9 scaled 950
\font=cmsl9 scaled 1050
\ifplaintex
\headline{\vbox to 0pt{\vskip -4.5mm\line{\small\phead\ifnum
\count0=\startpage ISSN 1364-0380 (on line)
1465-3060 (printed) \hfill {\pnum\folio}\else\ifodd\count0\def\\{ }%
\ifx\theshorttitle\relax\thetitle\else\theshorttitle\fi\hfill{\pnum\folio}
\else\def\\{ and }{\pnum\folio}\hfill\ifx\theshortauthors\relax\theauthors
\else\theshortauthors\fi\fi\fi}\vss}}
\footline{\vbox to 0pt{\vglue 0mm\line{\small\pfoot\ifnum\count0=\startpage
\copyright\ \gtp\hfill\else
\gt, Volume \thevolumenumber\ (\thevolumeyear)\hfill\fi}\vss
}}
\else
\makeatletter
\def\@oddhead{{\small\lhead\ifnum\count0=\startpage ISSN 1364-0380 (on line)
1465-3060 (printed) \hfill {\lnum\number\count0}\else\ifodd\count0
\def\\{ }\ifx\theshorttitle\relax \thetitle \else\theshorttitle\fi\hfill
{\lnum\number\count0}\else\def\\{ and }{\lnum\number\count0}
\hfill\ifx\theshortauthors\relax 
\theauthors\else\theshortauthors\fi\fi\fi}}\def\@evenhead{\@oddhead}
\def\@oddfoot{\small\lfoot\ifnum\count0=\startpage\copyright\ \gtp\hfill\else
\gt, Volume \thevolumenumber\ (\thevolumeyear)\hfill\fi}
\def\@evenfoot{\@oddfoot}
\makeatother
\fi

\newwrite\gtoutfile
\long\gdef\makeheadfile{  
{\def\\{, }\def\s{ }
\immediate\openout\gtoutfile head.xxx
\immediate\write\gtoutfile{Proxy-for: \ifx\theasciiauthors\relax
\theauthors\else\theasciiauthors\fi\s<\ifx\theasciiemail\relax\theemail\else\theasciiemail\fi>}
\immediate\write\gtoutfile{\noexpand\\}
\immediate\write\gtoutfile{Authors: \ifx\theasciiauthors\relax
\theauthors\else\theasciiauthors\fi}
{\def\\{ }\immediate\write\gtoutfile{Title: \ifx\theasciititle\relax
\thetitle\else\theasciititle\fi}}
\immediate\write\gtoutfile{Subj-class: GT or SG or MG etc}
\immediate\write\gtoutfile{MSC-class: \theprimaryclass\ifx\thesecondaryclass\relax\else, \thesecondaryclass\fi}
\immediate\write\gtoutfile{Journal-ref: Geom. Topol. \thevolumenumber
(\thevolumeyear) \startpage-\finishpage}
\immediate\write\gtoutfile{Comments: Published by Geometry and Topology at}
\immediate\write\gtoutfile{\s\s http://www.maths.warwick.ac.uk/gt/GTVol\thevolumenumber/paper\thepapernumber.abs.html}
\immediate\write\gtoutfile{\noexpand\\}
\immediate\write\gtoutfile{}
\ifx\theasciiabstract\relax
\immediate\write\gtoutfile{\theabstract}\else
\immediate\write\gtoutfile{\theasciiabstract}\fi
\immediate\write\gtoutfile{}
\immediate\write\gtoutfile{\noexpand\\}
\immediate\write\gtoutfile{}
\immediate\closeout\gtoutfile}}  

\def\maketitlepage{\maketitlep\makeheadfile}
\let\maketitle\maketitlepage

\lognumber{467}
\received{15 June 2004}
\volumenumber{8}\papernumber{37}\volumeyear{2004}
\pagenumbers{1361}{1384}
\revised{25 October 2004}
\published{25 October 2004}
\accepted{25 October 2004}
\proposed{Walter Neumann}
\seconded{Shigeyuki Morita, Joan Birman}

\usepackage{amsmath,amssymb, epsf}

\newcommand{\nc}{\newcommand}
\nc{\ens}{\ensuremath}
\nc{\dmo}{\DeclareMathOperator}
\nc{\nt}{\newtheorem}

\nt{thm}{Theorem}[section]
\nt{mthm}{Main Theorem}
\nt{cor}[thm]{Corollary}
\nt{minithm}[thm]{Proposition}
\nt{lem}[thm]{Lemma}
\nt{step}[thm]{Step}
\nt{fact}[thm]{Fact}
\nt{q}{Question}
\nt{rem}{Remark}
\nc{\prop}{Proposition}

\renewcommand{\rk}{{\rm rk}}
\nc{\cent}{\ens{Z}}
\nc{\kg}{\mathcal K}
\nc{\T}{\mathcal I}
\nc{\ps}{\phi_\star}
\nc{\pair}{{\mathcal P}}
\nc{\curv}{{\mathcal C}}
\dmo{\comm}{Comm}
\nc{\commk}{\comm(\kg)}
\nc{\fin}{G}
\nc{\bc}{\mathcal{B}}
\nc{\pst}{\widehat{\phi}_{\star}}
\nc{\uc}{\al\text{-}\be\text{-}\ga}
\nc{\ucp}{\uc'}
\nc{\ucpp}{\al'\text{-}\be\text{-}\ga'}
\nc{\ucppp}{\al'\text{-}\be'\text{-}\ga'}
\nc{\newc}{$\al'$-$\be$-$\ga'$ \hskip .01in}
\nc{\tri}{$\be^-\al\de$ \hskip .01in}
\nc{\bad}{$\be^\pm\al\de$ \hskip .01in}
\nc{\bae}{$\be^\pm\al\ep$ \hskip .01in}
\nc{\bcd}{$\be^\pm\ga\de$ \hskip .01in}
\nc{\bce}{$\be^\pm\ga\ep$ \hskip .01in}
\nc{\adce}{$\al\de\ga\ep$ \hskip .01in}
\nc{\bdce}{$\be^\pm\de\ga\ep$ \hskip .01in}
\nc{\bcea}{$\be^\pm\ga\ep\al$ \hskip .01in}
\nc{\bead}{$\be^\pm\ep\al\de$ \hskip .01in}
\nc{\badc}{$\be^\pm\al\de\ga$ \hskip .01in}
\nc{\badce}{$\be^\pm\al\de\ga\ep$ \hskip .01in}
\nc{\bcead}{$\be^\pm\ga\al\ep\de$ \hskip .01in}
\nc{\bcdae}{$\be^\pm\ga\de\al\ep$ \hskip .01in}
\nc{\baecd}{$\be^\pm\al\ep\ga\de$ \hskip .01in}
\nc{\bottom}{$\be^-$}
\nc{\topp}{$\be^+$}
\nc{\td}{\text{-}}

\nc{\al}{\ens{\alpha}}
\nc{\be}{\ens{\beta}}
\nc{\ga}{\ens{\gamma}}
\nc{\de}{\ens{\delta}}
\nc{\ep}{\ens{\epsilon}}
\nc{\Ga}{\ens{\Gamma}}
\nc{\ta}{\ens{\tilde{a}}}
\nc{\tb}{\ens{\tilde{b}}}
\nc{\De}{\ens{\Delta}}

\dmo{\mcg}{Mod} 
\dmo{\pmcg}{PMod} 
\nc{\mods}{\ens{\mcg(S)}}
\dmo{\diff}{Diff}
\dmo{\homeop}{\ens{Homeo^+}}
\dmo{\homeopm}{\ens{Homeo^\pm}}
\dmo{\homeo}{\ens{Homeo}}
\dmo{\gin}{i} 

\dmo{\cc}{C}
\nc{\ccs}{\cc(S)}
\nc{\ccos}{\cc^0(S)}
\nc{\csep}{\cc_s(S)}
\nc{\csepo}{\cc_s^0(S)}
\nc{\tg}{\mathcal T(S)}
\nc{\tgo}{\mathcal T^{0}(S)}
\nc{\tgeom}{\mathcal TG(S)}
\dmo{\si}{SI}

\dmo{\out}{Out}
\dmo{\aut}{Aut}
\nc{\zt}{\ens{{\mathbb Z}_2}}
\nc{\pslz}{\ensuremath{PSL_2(\z)}}
\nc{\spz}{\ens{Sp(2g,\z)}}
\nc{\z}{\ens{{\mathbb Z}}}
\nc{\br}{\ensuremath{{\mathbb R}}}
\dmo{\isom}{Isom} 

\nc{\lra}{\ens{\longrightarrow}}
\nc{\ra}{\ens{\rightarrow}}
\nc{\set}[1]{\{#1\}}
\nc{\genby}[1]{\langle #1 \rangle}

\nc{\bs}{\bigskip}
\nc{\bpf}{\begin{proof}}
\nc{\epf}{\end{proof}}
\nc{\p}[1]{\medskip{\bf #1}\qua\ignorespaces}

\nc{\bcom}{}

\nc{\bl}{ \begin{list}{$\cdot$}{
}
}

\nc{\el}{\end{list}}

\nc{\pic}[2]{\begin{figure}[htb] \center{\leavevmode 
\epsfbox{#1.eps}} \caption{#2} \label{#1pic}\end{figure}} 

\nc{\pics}[3]{\epsfysize=#3 cm \begin{figure}[htb] 
\center{\leavevmode \epsfbox{#1.eps}} \caption{#2}  
\label{#1pic}\end{figure}}

\begin{document}

\title{Commensurations of the Johnson kernel}

\author{Tara E Brendle\\Dan Margalit} 

\address{Department of Mathematics, Cornell University\\310 Malott Hall, 
Ithaca, NY  14853, USA}
\secondaddress{Department of Mathematics, University of Utah\\155 S 1440 East, Salt Lake City, UT 84112, USA}

\asciiaddress{Department of Mathematics, Cornell University\\310 Malott
Hall, Ithaca, NY 14853, USA\\and\\Department of Mathematics, University of
Utah\\155 S 1440 East, Salt Lake City, UT 84112, USA}

\asciiemail{brendle@math.cornell.edu, margalit@math.utah.edu}

\gtemail{\mailto{brendle@math.cornell.edu}, 
\mailto{margalit@math.utah.edu}}



\begin{abstract}
Let $\kg$ be the subgroup of the extended mapping class group,
$\mods$, generated by Dehn twists about separating curves.  Assuming
that $S$ is a closed, orientable surface of genus at least 4, we
confirm a conjecture of Farb that $\comm(\kg) \cong \aut(\kg) \cong
\mods$.  More generally, we show that any injection of a finite index
subgroup of $\kg$ into the Torelli group $\T$ of $S$ is induced by a
homeomorphism.  In particular, this proves that $\kg$ is co-Hopfian
and is characteristic in $\T$.  Further, we recover the result of Farb
and Ivanov that any injection of a finite index subgroup of $\T$ into
$\T$ is induced by a homeomorphism.  Our method is to reformulate
these group theoretic statements in terms of maps of curve
complexes.
\end{abstract}

\asciiabstract{%
Let K be the subgroup of the extended mapping class group, Mod(S),
generated by Dehn twists about separating curves.  Assuming that S is
a closed, orientable surface of genus at least 4, we confirm a
conjecture of Farb that Comm(K), Aut(K) and Mod(S) are all isomorphic.
More generally, we show that any injection of a finite index subgroup
of K into the Torelli group I of S is induced by a
homeomorphism.  In particular, this proves that K is co-Hopfian
and is characteristic in I.  Further, we recover the result of Farb
and Ivanov that any injection of a finite index subgroup of I into
I is induced by a homeomorphism.  Our method is to reformulate
these group theoretic statements in terms of maps of curve complexes.}

\maketitlepage

\section{Introduction}

The {\em extended mapping class group} of a surface $S$ is: \[ \mods =
\pi_0(\homeopm(S)) \] In this paper we show that certain topologically
defined subgroups of $\mods$ have algebraic structures which strongly
reflect their topological origins.  We assume throughout that $S$ is a
closed, oriented surface.  Also, unless specifically stated otherwise,
the genus $g$ of $S$ is assumed to be at least 4.

The {\em Torelli group} $\T=\T(S)$ is the subgroup of $\mods$
consisting of elements which act trivially on $H = H_1(S,\z)$, and
$\kg=\kg(S)$ denotes the subgroup of $\T$ generated by Dehn twists
about separating curves.  Birman and Powell initially suggested that
$\kg$ has finite index in $\T$.  However, Johnson, who first singled
out $\kg$ for study in its own right, showed (among other things) that
$\kg$ arises naturally as the kernel of a map from $\T$ to
$\wedge^3(H)/H$ \cite{djII}; hence we call $\kg$ the {\em Johnson
kernel}.

Building on Johnson's work, Morita showed that every homology 3--sphere
is obtained by splitting $S^3$ along a Heegaard surface and regluing
by an element of the Johnson kernel \cite{sm}.  Morita further
demonstrated that the natural map from $\kg$ to ${\mathbb Z}$, given
by computing the Casson invariant of the associated homology 3--sphere,
is a homomorphism (this is not true for the analagous map from $\T$ to
${\mathbb Z}$).  More recently, Biss and Farb applied ideas of
McCullough and Miller \cite{mm} to prove that $\kg$ is not finitely
generated \cite{bifa}.

\p{Abstract commensurators} The {\em abstract commensurator} of a
group $\Ga$, denoted $\comm(\Ga)$, is the group of isomorphisms of
finite index subgroups of $\Ga$ (under composition), with two such
isomorphisms equivalent if they agree on a finite index subgroup of
$\Ga$.  The product of $\phi \co  G \to H$ with $\psi \co  G' \to H'$ is a
map defined on $\phi^{-1}(H \cap G')$.

While there is always a map from $\aut(\Ga)$ to $\comm(\Ga)$, one
generally expects $\comm(\Ga)$ to be much larger (for instance, if
$\Ga$ is arithmetic).  As an example, $\comm({\mathbb Z}^n) \cong
\textrm{GL}_n({\mathbb Q})$, whereas $\aut({\mathbb Z}^n) \cong
\textrm{GL}_n({\mathbb Z})$.  On the other hand, Ivanov used
Theorem~\ref{autcs} below to prove $\comm(\mods) \cong \mods$
\cite{ni}, Farb and Mosher showed that any word-hyperbolic
surface-by-free group has finite index in its abstract commensurator
\cite{fm}, and Farb and Handel proved $\comm(\out(F_n)) \cong
\out(F_n)$ \cite{fh}.

We confirm a conjecture of Farb, from his talk at an AMS sectional
meeting in 2002 (Question 2 of \cite{bf}):

\begin{mthm}\label{comm}For $S$ a surface of genus $g \geq 4$, we have: \[ \comm(\kg) \cong \aut(\kg) \cong \mods \] \end{mthm}

This follows from a more general theorem:

\begin{mthm}\label{intorelli}Let $S$ be a surface of genus $g \geq 4$.  Any injection $\phi\co  \fin \to \T$, where $\fin$ is a finite index subgroup of $\kg$, is induced by an element $f$ of $\mods$ in the sense that $\phi(h) = fhf^{-1}$ for all $h \in \fin$.\end{mthm}

\p{Co-Hopfian property} A group is {\em co-Hopfian} if each of its
injective endomorphisms is an isomorphism.  In general, this is more
rare (and harder to prove) than the Hopfian property, which, for
instance, holds for all finitely generated linear groups by work of
Mal'cev (recall that a group is Hopfian if every surjective
endomorphism is an isomorphism).

The co-Hopfian property has been demonstrated for lattices in
semisimple Lie groups by Prasad \cite{gp}, torsion-free non-elementary
freely indecomposable word-hyperbolic groups by Sela \cite{zs},
$\mods$ by Ivanov and McCarthy \cite{im}, the braid group modulo its
center by Bell and Margalit \cite{bm}, and finite-index subgroups of
$\out(F_n)$ by Farb and Handel \cite{fh} (the case of $\out(F_n)$
itself follows from work of Bridson and Vogtmann \cite{bv}).

As consequences of Main Theorem~\ref{intorelli}, we establish two
basic algebraic properties of $\kg$.  First, $\kg$ is co-Hopfian.
Also, $\kg$ is characteristic in $\T$.  The former was asked by Farb
as Question 3 in his lecture \cite{bf}.

\begin{cor}\label{cohopf}If $S$ is a surface of genus $g \geq 4$, then $\kg$ and all of its finite index subgroups are co-Hopfian.\end{cor}

To obtain the corollary for a finite index subgroup $\fin$, note that
given an injection $\fin \to \fin$, the automorphism of $\kg$ supplied
by Main Theorem~\ref{intorelli} sends each of the (finitely many)
cosets of $\fin$ in $\kg$ entirely into another such coset.

\begin{cor}\label{cistic}If $S$ is a surface of genus $g \geq 4$, then $\kg$ is characteristic in $\T$.  Finite index subgroups of $\kg$ are characteristic in $\T$ up to conjugacy.\end{cor}

By ``characteristic up to conjugacy'', we mean that automorphisms of
$\T$ preserve the conjugacy class in $\T$ of a finite index subgroup
of $\kg$.

Our methods, which are largely inspired by the work of Farb and
Ivanov, can be used recover their main results:

\begin{thm}\label{auttor}For $S$ a surface of genus $g \geq 4$, we have: \[ \comm(\T) \cong \aut(\T) \cong \mods \] \end{thm}

\begin{thm}\label{cohopftor}If $S$ is a surface of genus $g \geq 4$, then $\T$ and all of its finite index subgroups are co-Hopfian.  Further, any injection $\psi\co \fin \to \T$, where $\fin$ is a finite index subgroup of $\T$, is induced by an element $f$ of $\mods$ in the sense that $\psi(h)=fhf^{-1}$ for all $h \in \fin$.\end{thm}

We remark that Farb and Ivanov proved these for $g \geq 5$, and
McCarthy and Vautaw extended the result that $\aut(\T) \cong \mods$ to
$g \geq 3$ \cite{mv}.  We also note that in the case of the Johnson
kernel itself, Corollary~\ref{cistic} also follows from
Theorem~\ref{cohopftor}, since conjugation by an element of $\mods$
preserves the set of Dehn twists about separating curves (see
Fact~\ref{conjtwist} in Section~\ref{bg}).

Our basic strategy is to recast algebraic statements about $\kg$ in
terms of combinatorial topology.  That is, given an injection of a
finite index subgroup of $\kg$ into $\T$, we construct a map of the
complex of separating curves $\csep$, which has vertices corresponding
to separating curves in $S$ and edges corresponding to disjointness
(see Section~\ref{bg}).  We then prove:

\begin{thm}\label{autsep}For $S$ a surface of genus $g \geq 4$, we have: \[ \aut(\csep) \cong \mods \] \end{thm}

More generally, we look at superinjective maps of $\csep$.
Superinjective maps were defined and used by Irmak in order to show
that, for most surfaces, any injection of a finite index subgroup of
$\mods$ into $\mods$ is induced by a homeomorphism \cite{ei}
\cite{eiII}.  See Section~\ref{bg} for a precise definition.

\begin{thm}\label{superthm}If $S$ is a surface of genus $g \geq 4$, then every superinjective map of $\csep$ is induced by an element of $\mods$.\end{thm}

The corresponding results for Torelli groups concern the Torelli
complex $\tg$, which has a vertex for each separating curve and
bounding pair in $S$, and edges between vertices which can be realized
disjointly in $S$ (see Section~\ref{bg}).  This complex was originally
defined by Farb and Ivanov (their Torelli geometry $\tgeom$ is $\tg$
with a certain marking) \cite{fi}.

We prove:

\begin{thm}\label{auttg}For $S$ a surface of genus $g \geq 4$, we have: \[ \aut(\tg) \cong \mods \] \end{thm}

\begin{thm}\label{supertg}If $S$ is a surface of genus $g \geq 4$, then any superinjective map of $\tg$ is induced by an element of $\mods$.\end{thm}

The previous two theorems generalize the following theorem of Farb and
Ivanov about automorphisms of $\tgeom$ (here, automorphisms preserve
the marking) \cite{fi}:

\begin{thm}For $S$ a surface of genus $g \geq 5$, we have: \[ \aut(\tgeom) \cong \mods \]\end{thm}

All of the above results concerning maps of ``curve complexes'' are
based on the seminal theorem of Ivanov about the complex of curves
$\ccs$ (defined in Section~\ref{bg}) \cite{ni}:

\begin{thm}\label{autcs} For $S$ a surface of genus $g \geq 3$, we have: \[ \aut(\ccs) \cong \mods \]\end{thm}

In lower genus cases, explored by Korkmaz and Luo \cite{mk} \cite{fl},
there are exceptions to Ivanov's theorem.

\p{Remarks} Despite the deep connection between mapping class groups
and arithmetic groups, none of the groups $\mods$, $\T$, and $\kg$ are
arithmetic \cite{fi} \cite{ni}.  For $\kg$, it suffices to note that
arithmetic groups are finitely generated and compare with the
aforementioned theorem of Biss and Farb.

Mustafa Korkmaz has pointed out another approach to computing
$\aut(\csep)$, using his theorem that $\aut(\ccs) \cong \mods$ for
punctured spheres.  However, since there is currently no analog of
Irmak's superinjectivity result for genus zero, this method does not
seem to recover our second main theorem.

\p{Outline} We start with any injection: \[ \phi \co  \fin
\hookrightarrow \T \] where $\fin$ is any finite index subgroup of
$\kg$.  As per Main Theorem~\ref{intorelli}, we aim to produce a
mapping class $f$ which induces $\phi$.

{\bf Step 1}\qua {\sl The injection $\phi$ induces a map $\ps$ of the vertices
of $\csep$ into $\tg$.}  We show that $\phi$ takes a Dehn twist to
either a Dehn twist or a ``bounding pair map'' (a product of two Dehn
twists), so we get a natural map $\ps$ from curves to collections of
curves.  This follows from Proposition~\ref{twistpres}, which relies
on Ivanov's algebraic characterization of ``simple'' mapping
classes---the centers of their centralizers are infinite cyclic.

{\bf Step 2}\qua {\sl The set map $\ps$ extends to a simplicial map from
$\csep$ to $\tg$.}  This follows from Lemma~\ref{disjlem}, which uses
the fact that Dehn twists commute if and only if the corresponding
curves are disjoint.

{\bf Step 3}\qua {\sl The map $\ps$ is really a map of $\csep$ to itself.}  In
other words, we show in Proposition~\ref{genuslem} that $\ps$ takes
separating curves to separating curves.  The idea is that the rank of
a maximal abelian subgroup on either side of a separating curve is
even, and on either side of a bounding pair is odd.  We use the
assumption that $S$ is closed in Section~\ref{twistsec}.

{\bf Step 4}\qua {\sl The map $\ps$ extends to a map $\pst$ of $\ccs$.}  In
Sections~\ref{spsec} and~\ref{3c}, we define an action on
nonseparating curves.  The key is that a nonseparating curve is
characterized by the separating curves which cut off a genus 1
subsurface containing that curve.  To define $\pst$ on a particular
nonseparating curve, it suffices to choose two such separating
curves---a ``sharing pair''.  We use the assumption on genus here.

To show that $\pst$ is well-defined, ie, that the map is independent
of choices of sharing pairs, we apply the chain-connectedness of a
certain arc complex defined by Harer.

{\bf Step 5}\qua {\sl The map $\pst$ is a superinjective simplicial map of
$\ccs$ to itself.}  As we shall see in Section~\ref{homeo}, this
follows easily from the superinjectivity of $\ps$ and the definition
of $\pst$.  It is then a theorem of Irmak that $\pst$ is induced by
some $f \in \mods$.

{\bf Step 6}\qua {\sl Both $\ps$ and $\phi$ itself are induced by $f$.}  We
accomplish this by showing that $\phi(h)$ and $fhf^{-1}$ have the same
action on $\csep$ for any $h \in G$.  Thus, we conclude in
Section~\ref{mainproof} that $\phi$ is really the restriction of an
element of $\aut(\kg)$ and $\ps$ is really an automorphism of $\csep$.
To summarize:
\[ \phi \mapsto \ps \mapsto \pst \mapsto f \]
At this point, we have established Main Theorem~\ref{intorelli} and
Theorem~\ref{superthm}, and we have: \[ \aut(\kg) \cong \comm(\kg) \to
\mods \]

The natural map $\eta\co  \mods \to \aut(\kg)$ is an inverse, so it is an
isomorphism.

In Section~\ref{reproof}, we explain how to get the related results
about Torelli groups.

\p{Acknowledgements} The authors are indebted to Bob Bell, Mladen
Bestvina, Nate Broaddus, Indira Chatterji, Benson Farb, Karen
Vogtmann, and Kevin Wortman for reading our drafts, and for useful
conversations.  In particular, Benson Farb suggested that we extend
our earlier results to finite index subgroups.  Bob Bell was
especially helpful in listening to many of our proofs.  We are
grateful to Mustafa Korkmaz for sharing his alternate approach to the
problem, and for helpful suggestions.  We are also very thankful for
corrections made by Walter Neumann to an earlier manuscript.  Finally,
we thank Martin Kassabov for pointing out that our proof works in the
added case of genus four.

Both authors would like to thank Joan Birman, Benson Farb, and Kevin
Wortman for their enthusiasm and encouragement on this project.  We
would also like to express our gratitude to Cornell University and the
University of Utah for pleasant working environments during our
respective visits.  Finally, we are grateful to the referee for
carefully reading the paper and suggesting many improvements.

The first author is partially supported by a VIGRE postdoc under NSF
grant number 9983660 to Cornell University.  The second author is
partially supported by a VIGRE postdoc under NSF grant number 0091675
to the University of Utah.

\section{Background}\label{bg}

\subsection{Curves}When no confusion arises, we will use {\em curve} to mean either ``simple closed curve'' or ``isotopy class of simple closed curves''.

A curve $c$ is {\em separating} if $S-c$ is not connected.  The {\em
genus} of a separating curve $c$ on $S$ is the minimum of the genera
of the two components of $S-c$.

\p{Bounding pairs} A {\em bounding pair}, denoted $(a,b)$, is an
ordered pair of disjoint, homologous nonseparating curves.

\p{Intersection}We use $\gin(\cdot,\cdot)$ to mean geometric
intersection number.  If either term is a bounding pair, we use the
following formula: \[ \gin((a,b),\cdot)=\gin(a,\cdot)+\gin(b,\cdot) \]
We say that two curves and/or bounding pairs {\em intersect} if their
geometric intersection is nonzero.

\subsection{Generators}

A {\em Dehn twist} about a simple closed curve $c$, denoted $T_c$, is
the isotopy class of a homeomorphism which is the identity outside of
a regular neighborhood $N$ of $c$, and which is described on the
annulus $N$ by Figure~\ref{3dtwistpic}.

\pics{3dtwist}{A Dehn twist}{3}

\p{Bounding pair maps}Johnson proved that $\T$ is generated by
finitely many {\em bounding pair maps} \cite{djI}, which are given by:
\[ BP(a,b) = T_aT_b^{-1} \] for $(a,b)$ a bounding pair.  The Johnson
kernel contains no bounding pair maps.

\p{Rank} Both Dehn twists and bounding pair maps can appear in free
abelian subgroups of maximal rank in $\T$.  In general, we denote by
$\rk(\Ga)$ the maximal rank of a free abelian subgroup in a group
$\Ga$.  The next proposition follows from the classification of
abelian subgroups of $\mods$ obtained by Birman, Lubotzky, and
McCarthy \cite{blm}.

\begin{minithm}\label{max}For any surface $S$, $\rk \kg = \rk \T$.  If $S$ has genus $g$ and $b$ boundary components where $b \in \set{0,1,2}$, then $\rk \kg = \rk \T = 2g-3+b$.\end{minithm}

In particular, Proposition~\ref{max} tells you the size of a maximal
collection of disjoint separating curves on $S$, as there is always a
subgroup of rank $\rk \kg = \rk \T$ generated by Dehn twists, and any
collection of $n$ disjoint curves gives rise to a free abelian
subgroup of rank $n$.  The parities of the ranks in
Proposition~\ref{max} are crucial in Proposition~\ref{genuslem}.

\pics{max}{A maximal collection of disjoint separating curves}{6}

\subsection{Relations}

While a complete set of relations for $\kg$ is not known, the results
of this paper only rely on the commuting relation of
Fact~\ref{commute} below.

\begin{fact}\label{commutecrs}Let $f, T_c \in \T$, where $T_c$ is a Dehn twist.  For $k \neq 0$, we have $[f,T_c^k]=1$ if and only if $f(c)=c$.\end{fact}

The last statement follows from the fact that elements of $\T$ are
orientation preserving, together with the identity:

\begin{fact}\label{conjtwist}For $f, T_c \in \mods$, we have $fT_c^kf^{-1}=T_{f(c)}^{\ep k}$, where $\ep$ is $1$ if $f$ is orientation preserving, and $-1$ otherwise.\end{fact}

As a special case of the Fact~\ref{commutecrs}, we have:

\begin{fact}\label{commute}For $f \in \set{BP(a,b),T_c}$, $h \in \set{BP(x,y),T_z}$, and $j$ and $k$ nonzero, we have $[f^j,h^k]=1$ if and only if the intersection number between the corresponding curves and/or bounding pairs is zero.\end{fact}

\subsection{Curve complexes} We will require the use of three different abstract complexes which have vertices corresponding to curves (or bounding pairs), and have edges between vertices which can be realized disjointly.

\p{Complex of curves} The {\em complex of curves} $\ccs$ for a surface
$S$, defined by Harvey \cite{wjh}, is the abstract simplicial flag
complex with a vertex for each curve in $S$ and edges between vertices
which can be realized as disjoint curves in $S$.

\p{Complex of separating curves} The subcomplex of $\ccs$ spanned by
vertices corresponding to separating curves is called the {\em complex
of separating curves}, and is denoted $\csep$.

Farb and Ivanov introduced $\csep$ and showed that it is connected for
all surfaces of genus at least 3 \cite{fi}, in particular for the
surfaces we consider (for alternate proofs, see \cite{ms} and
\cite{mv}).

\p{Torelli complex}We denote by $\tg$ the abstract simplicial flag
complex with vertices corresponding to separating curves and bounding
pairs in $S$, and edges between vertices that can be realized
disjointly in $S$.  This {\em Torelli complex} is the simplicial
complex underlying the Torelli geometry of Farb and Ivanov \cite{fi}.
It follows from the connectedness of $\csep$ that $\tg$ is connected.

\p{Superinjective maps} A simplicial map $\psi$ of $\ccs$, $\csep$, or
$\tg$ into itself is {\em superinjective} if $\gin(\psi(a),\psi(b))
\neq 0$ whenever $\gin(a,b) \neq 0$.  Here, we are computing
intersection number on the curves (and/or bounding pairs)
corresponding to the vertices of the complex.

It is not hard to prove that superinjectivity implies injectivity for
all of the above complexes, but this does not follow immediately from
the definition (see \cite{ei}).  It does follow from the definition
that superinjective maps preserve disjointness (the maps are
simplicial).

\p{Remark} If $X$ is a flag complex, then $\aut(X) = \aut(X^{(1)})$,
where $X^{(1)}$ is the 1--skeleton of $X$.  Thus, we may restrict our
attention to vertices and edges when discussing automorphisms of
$\ccs$, $\csep$, and $\tg$.

\section{Separating curves}\label{twistsec}

Let $\phi\co  \fin \to \T$ be an injection, where $\fin$ is a finite
index subgroup of $\kg$.  The goal of this section is to prove that
the image under $\phi$ of a power of a Dehn twist in $\fin$ is a power
of a Dehn twist.  This will imply two things: \begin{enumerate}
\item $\phi(\fin) < \kg$
\item $\phi$ induces a superinjective map $\ps$ of $\csep$.
\end{enumerate}
We first show that $\phi$ takes a power of a Dehn twist to a power of
either a Dehn twist or a bounding pair map.  We then rule out bounding
pair maps.

\subsection{Centers of centralizers}

We now give a series of three propositions: an algebraic
characterization of Dehn twists in $\fin$ in terms of centers of
centralizers, a comparison between the centers of centralizers of an
element and its image under an injection, and a characterization of
Dehn twists (and bounding pair maps) in $\T$.

The first proposition follows from work of Ivanov which in turn rests
on Thurs\-ton's theory of surface homeomorphisms \cite{ninotes}.  First,
we will require a lemma, which follows immediately from the fact that
if a mapping class permutes components of a surface, it necessarily
acts nontrivially on the first homology of $S$.

\begin{lem}\label{pure} If $f \in \T$, and $f(c)=c$ for some curve $c$, then $f$ leaves invariant the components of $S-c$.\end{lem}

\begin{minithm}\label{rankkg}If $f \in \fin$ is a power of a Dehn twist, then $\cent(C_\fin(f)) \cong \z$ and $\rk C_\fin(f) = 2g-3$.\end{minithm}

\bpf

Let $f = T_c^k \in \fin$ be a power of a Dehn twist. First, we know
that $\rk C_\fin(f) = 2g-3$ by Proposition~\ref{max}.  It remains to
show that $\cent(C_\fin(f)) \cong \z$.

Each element of $C_\fin(f)$ fixes the separating curve $c$
(Fact~\ref{commutecrs}), and has well-defined restrictions to $S-c$ by
Lemma~\ref{pure}. Suppose $h \in C_\fin(f)$ has nontrivial support on
a component $S'$ of $S-c$ (note that such a component necessarily has
genus greater than one).  Then there is a separating curve $a$ in $S'$
with $h(a) \neq a$.  By Fact~\ref{commutecrs}, $h$ does not commute
with nontrivial powers of $T_a$ (some of which are in $C_\fin(f)$),
and hence $h \notin \cent(C_\fin(f))$.  It follows that
$\cent(C_\fin(f))$ can only consist of powers of $T_c$.
\epf

The following proposition is due to Ivanov and McCarthy \cite{im}:

\begin{minithm}\label{gpthy}
Let $\phi\co  \fin \to \fin'$ be an injection between groups with $\rk
(\fin)$ $= \rk (\fin') < \infty$.  Let $H \cong \z^{\rk(\fin)}$ be a
subgroup of $\fin$, and let $h \in H$.  Then
\[
\rk Z(C_{\fin'}(\phi(h)))\leq \rk Z(C_\fin(h))
\]
\end{minithm}

We now require the following proposition of Farb and Ivanov (private
communication; announced in \cite{fi}).

\begin{minithm}\label{trank}An element $h$ of $\T$ is a power of a Dehn twist or a bounding pair map if and only if both of the following hold: 
\begin{enumerate} 
\item $\cent(C_\T(h)) \cong \z$ 
\item $\rk C_\T(h) = 2g-3$
\end{enumerate}
\end{minithm}

We now put the above results together:

\begin{minithm}\label{twistpres}If $\phi\co  \fin \to \T$ is an
injection for $\fin$ a finite index subgroup of $\kg$, and $f=T_c^k
\in \fin$, then $\phi(f)$ is a power of either a Dehn twist or a
bounding pair map.\end{minithm}

\bpf

We will show that $\phi(f)$ satisfies the two conditions of
Proposition~\ref{trank}.  We have that $\rk \cent(C_\T(\phi(f))) \leq
\rk \cent(C_\kg(f)) = 1$, by Propositions~\ref{max}, \ref{rankkg},
and~\ref{gpthy}.  Since the center of the centralizer of an element
$h$ in a group always contains $h^n$ for all integers $n$, and since
$\T$ is torsion free, we have that $\rk \cent(C_\T(\phi(f)))$ is
precisely $1$.

Since $\phi$ is an injection, we know that $\rk C_\T(\phi(f)) \geq \rk
C_\kg(f) = 2g - 3$, with the latter equality given by
Proposition~\ref{rankkg}.  Proposition~\ref{max} now implies that this
inequality must in fact be an equality, which proves the proposition.
\epf

By \prop~\ref{twistpres}, $\phi$ gives a natural map $\ps$ from the
vertices of $\csep$ to the vertices of $\tg$ defined by: \[
\phi(T_a^k)=T_{\ps(a)}^{k'} \qquad \text{ or } \qquad \phi(T_a^k) =
BP(\ps(a))^{k'} \] In the remainder of this section, we will argue
that the latter is an impossibility, and further that $\ps$ is a
superinjective map of $\csep$.

\subsection{Basic topology}\label{basic}

We will now show that $\ps$ preserves some elementary topological
relationships between separating curves on $S$.  In particular, it
will follow from Lemma~\ref{disjlem} and Proposition~\ref{genuslem}
that $\ps$ is a superinjective map of $\csep$.

\begin{lem}[Disjointness]\label{disjlem}If $a$ and $b$ are separating curves in $S$, then $\gin(a,b) \neq 0$ if and only if $\gin(\ps(a),\ps(b)) \neq 0$.\end{lem}

The above lemma follows immediately from Fact~\ref{commute} and the
assumption that $\phi$ is an injective homomorphism.

We define a {\em side} of a separating curve $z$ (or bounding pair
$(a,b)$) to be a component of $S-z$ (or $S-(a \cup b)$).

\begin{lem}[Sides]\label{sides}If separating curves $a$ and $b$ are on the same side of a separating curve $z$, then $\ps(a)$ and $\ps(b)$ are separating curves or bounding pairs on the same side of $\ps(z)$.\end{lem}
\bpf

If the curves $a$ and $b$ are on the same side of $z$, then there is a
curve $c$ which intersects $a$ and $b$ but not $z$.  Apply
Lemma~\ref{disjlem}.
\epf

\begin{minithm}[Genus]\label{genuslem} If $z$ is a genus $m$ separating curve, then $\ps(z)$ is a genus $m$ separating curve.  Further, if $a$ is on a genus $m$ side of $z$, then $\ps(a)$ is on a genus $m$ side of $\ps(z)$.\end{minithm}

\bpf

First, if $z$ is a genus $m$ curve, then any maximal collection of
disjoint separating curves in $S$ which contains $z$ is of the form:
\[ \set{a_1,\dots,a_{2m-2},z,b_1, \dots, b_{2(g-m)-2}} \] where the
$a_i$ are on one side of $z$, and the $b_i$ are on the other side of
$z$ (see Proposition~\ref{max}).

By Proposition~\ref{twistpres}, Lemma~\ref{disjlem}, and
Proposition~\ref{max}, we have that the set
$\set{\ps(a_i),\ps(z),\ps(b_i)}$ is a maximal collection of mutually
disjoint separating curves and bounding pairs.  By Lemma~\ref{sides}
and Proposition~\ref{max}, we either have that $\ps(z)$ is a genus 1
separating curve and the $\ps(a_i)$ and $\ps(b_i)$ are on the same
side of $\ps(z)$, or $\ps(z)$ is a genus $m$ separating curve and the
$\ps(a_i)$ and $\ps(b_i)$ are on different sides of $\ps(z)$.

If $m=1$, both cases are the same, so the image under $\ps$ of $z$ (or
any genus 1 curve) must be a genus 1 curve.  Since there are always
exactly $g$ genus 1 curves in any maximal collection of separating
curves and bounding pairs, we can only have that $\ps(z)$ is a genus 1
curve if $z$ was a genus 1 curve to begin with.  Both statements of
the proposition follow.
\epf

Note that the above proof does not work in general for surfaces with
boundary.

\p{Action on separating curves} Applying Proposition~\ref{genuslem}
and Lemma~\ref{disjlem}, we now have that $\ps$ is a superinjective
map of $\csep$ defined by: \[ \phi(T_c^k) = T_{\ps(c)}^{k'} \]

\section{Nonseparating curves}\label{nonsepsec}

We show in this section that any superinjective map $\ps$ of $\csep$
(in particular, the one defined in Section~\ref{twistsec}) can be
extended to a map $\pst$ of $\ccs$; that is, we get a natural map on
nonseparating curves as well.  We will then see in Section~\ref{homeo}
that $\pst$ in turn gives rise to an element of $\mods$.

A nonseparating curve $\be$ is uniquely determined by a pair of genus
1 separating curves with the property that $\be$ lies on both of the
corresponding genus 1 subsurfaces.  This notion gives rise to the
following definition.

\p{Sharing pairs} Let $a$ and $b$ be genus 1 curves bounding genus 1
subsurfaces $S_a$ and $S_b$ of $S$.  We say that $a$ and $b$ {\em
share} a nonseparating curve $\be$ if $S_a \cap S_b$ is an annulus
containing $\be$ (see Figure~\ref{sppic}) and $S-(S_a \cup S_b)$ is
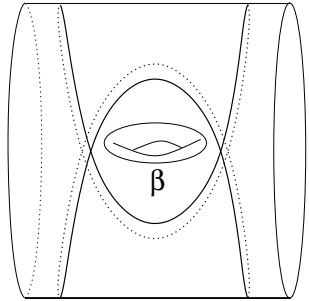
connected.  We also say that $a$ and $b$ form a {\em sharing pair} for
$\be$.

\pics{sp}{A sharing pair for $\be$}{4}

\p{Action on nonseparating curves} The map $\pst$ is defined on
nonseparating curves as follows.  If $\pair(\be)=\set{a,b}$ is a
sharing pair for a nonseparating curve $\be$, then $\pst(\be)$ is the
curve shared by $\ps(\pair(\be))$, which is defined as
$\set{\ps(a),\ps(b)}$.

To see that $\pst$ is well-defined on nonseparating curves, we must
check that $\ps(\pair(\be))$ is a sharing pair (Section~\ref{spsec}),
and that $\pst(\be)$ is independent of choice of $\pair(\be)$
(Section~\ref{3c}).

\subsection{Preserving sharing pairs}\label{spsec}

We now give a characterization of sharing pairs.  By the results of
Section~\ref{basic}, all the properties used in the characterization
are preserved by $\ps$.  The ideas here are inspired by related work
of Ivanov (see \cite{ni}, Lemma 1).

\begin{lem}[Sharing pairs]\label{nspchar}Let $a$ and $b$ be genus 1 curves in $S$.  Then $a$ and $b$ are a sharing pair if and only if there exist separating curves $w$, $x$, $y$, and $z$ in $S$ with the following properties: 
\bl
\item $z$ is a genus 2 curve bounding a genus 2 subsurface $S_z$.
\item $a$ and $b$ are in $S_z$ and intersect each other.
\item $x$ and $y$ are disjoint.
\item $w$ intersects $z$, but not $a$ and not $b$.
\item $x$ intersects $a$ and $z$, but not $b$.
\item $y$ intersects $b$ and $z$, but not $a$. 
\el
\end{lem}

\bpf

One direction is proven by construction: if $a$ and $b$ are a sharing
pair, then the other curves can be chosen as in Figure~\ref{nsppic}.

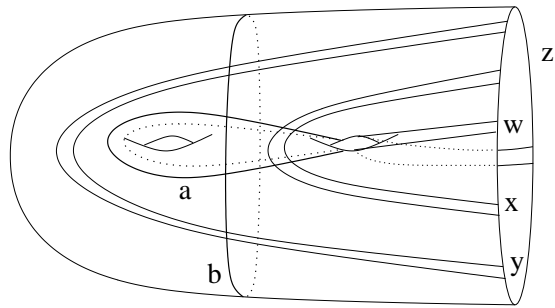
\pics{nsp}{Characterizing curves for a sharing pair}{4}

For the other direction, we restrict our attention to $S_z$.  On this
subsurface, $a$ and $b$ are genus 1 separating curves, and each of
$w$, $x$, and $y$ is a collection of arcs which separates $S_z$.  Let
$S_{az}$ denote the genus 1 surface of $S_z$ bounded by $a$ and $z$.

{\bf Step 1}\qua The arcs of $w$ are all parallel and nonseparating.

On $S_{az}$, we think of $b$ as a collection of disjoint arcs with
both endpoints on $a$, and of $w$ as a collection of disjoint arcs
with both endpoints on $z$.  By the hypotheses, the arcs of $b$ do not
intersect the arcs of $w$.

Now, these arcs of $b$ and $w$ are nontrivial in $H_1(S_{az},a)$ and
$H_1(S_{az},z)$, respectively: if any arcs of $w$, say, are all
parallel to $z$, then $w$ and $z$ are not in minimal position; if they
enclose a disc with one hole (namely, the hole is $a$), then there is
no room for the arcs of $b$.

If $w$ consists of non-parallel arcs, then $S_{az}-w$ is a disk with a
hole, so there is no room for the arcs of $b$.

{\bf Step 2}\qua $a$ and $b$ form a sharing pair.

Let $S_{wz}$ be the genus 1 subsurface of $S_z$ in the complement of
$w$.  Each of the two boundary components of $S_{wz}$ is made up of
exactly one arc from $z$ and one arc from $w$.  Since $a$ and $b$ are
genus 1 curves, they each split $S_{wz}$ into a torus with a hole and
a disk with two holes (the boundary components of $S_{wz}$).

As above, at least some arcs of $x$ and $y$ must be nontrivial in the
relative homology of $S_{wz}$, for otherwise the curves are not in
minimal position, or they are all parallel to the arcs of $w$, which
contradicts the fact that there are curves ($a$ and $b$) which
intersect $x$ and $y$, but not $w$.  We can thus ignore the arcs of
$x$ and $y$ which are parallel to $w$.

We will presently use the fact that all genus 1 curves on $S_{wz}$ (in
particular, $a$ and $b$) are obtained by taking the boundary of a
regular neighborhood in $S_{wz}$ of the union of $z \cup w$ with a
{\em defining arc} connecting the boundary components of $S_{wz}$.

If the defining arcs for $a$ and $b$ intersect, then there are no
choices for collections of arcs $x$ and $y$ which do not intersect
each other.  Therefore, the defining arcs for $a$ and $b$ do not
intersect (there is only one way to achieve this, topologically), and
it follows that $a$ and $b$ are a sharing pair.
\epf

We remark that the collection of curves in the lemma requires genus at
least 4, hence the hypothesis in our theorems.

Lemma~\ref{nspchar} implies:

\begin{minithm}\label{sharethm}Suppose $\ps$ is a superinjective map of $\csep$.  If two genus 1 curves $a$ and $b$ in $S$ form a sharing pair, then so do $\ps(a)$ and $\ps(b)$.\end{minithm}

\bpf

Since $a$ and $b$ share a curve, there are characterizing curves $w$,
$x$, $y$, and $z$ as in Lemma~\ref{nspchar}.  Since each property of
this collection of curves (disjointness, sides, genus) is preserved
(Lemmas~\ref{disjlem}--\ref{sides} and Proposition~\ref{genuslem}),
Lemma~\ref{nspchar} implies that $\ps(a)$ and $\ps(b)$ share a curve.
\epf

\subsection{Well-definedness}\label{3c}

We now have a map from the set of superinjective maps of $\csep$ to
the set of maps of $\ccs$, given by $\ps \mapsto \pst$.  In order to
show that this map is well-defined with respect to choice of sharing
pairs, it will be more convenient for us to consider sharing pairs
indirectly, via their ``spines''.

\p{Spines}Given two nonseparating curves $\al$ and $\be$ with
$\gin(\al,\be)=1$, we define $\bc(\al,\be)$ to be the genus 1
separating curve which is the boundary of a regular neighborhood of
$\al \cup \be$.  An ordered collection of three distinct curves
$\set{\al,\be,\ga}$ forms a {\em spine} for a sharing pair $\set{a,b}$
if: \[ \gin(\al,\be) = \gin(\be,\ga) = 1 \] \[ \gin(\al,\ga) \leq 1 \]
\[ \bc(\al,\be) = a \qquad \bc(\be,\ga) = b \] \[ S - (\al \cup \be
\cup \ga)\text{ is connected} \]

We denote this spine by $\al\td\be\td\ga$.

We can always choose a spine for a given sharing pair $\set{a,b}$,
although this choice is not unique.  If $\uc$ is a spine for a sharing
pair with $\gin(\al,\ga) = 0$, then $\al\td\be\td T_\be(\ga)$ and
$\al\td\be\td T_\be^{-1}(\ga)$ are also spines for $\set{a,b}$.  We
observe that every other spine for $\set{a,b}$ is obtained by applying
powers of $T_\be$ (curve by curve) to these 3 spines.  However, we
will not use this fact.

\p{Moves}We define a {\em move} between spines to be a change of the
form: \[ \uc \mapsto \ucp \] where $\ga\td\be\td\ga'$ is also a spine.
We remark that only one curve in the corresponding sharing pairs is
changed.  Note that the curves $\bc(\al,\be)$, $\bc(\be,\ga)$, and
$\bc(\be,\ga')$ form three sharing pairs for $\be$, since $\uc$,
$\ucp$ and $\ga\td\be\td\ga'$ are all spines.

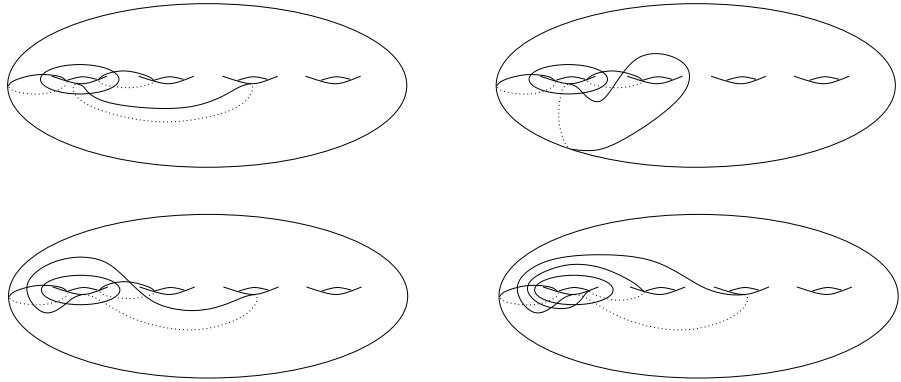
\pics{twistedmoves}{Moves on spines}{5}

A move is characterized by how many of $\gin(\al,\ga)$,
$\gin(\al,\ga')$, and $\gin(\ga,\ga')$ are 1, and how many are 0.
Using the nonseparating property of spines, it is straightforward to
prove:

\begin{lem}\label{4types}Any move is topologically equivalent to one of the 4 moves shown in Figure~\ref{twistedmovespic}.\end{lem}

The proof of the following lemma borrows from work of Ivanov and
McCarthy (Lemma 10.11 of~\cite{im}).  We will make repeated use of the
following fact: if $a$ is a genus 1 separating curve, and $\al$ and
$\be$ are any two distinct curves on the genus 1 surface bounded by
$a$, then a curve $z$ intersects $a$ if and only if it intersects at
least one of $\al$ and $\be$.

\begin{lem}\label{sharingpreserved}Let $\be$ be a nonseparating curve in $S$.  If $\set{a,b}$ and $\set{a,b'}$ are sharing pairs for $\be$ which have spines differing by a move, then $\set{\ps(a),\ps(b)}$ and $\set{\ps(a),\ps(b')}$ are sharing pairs for the same curve.\end{lem}

\bpf

Suppose that $\uc$ and $\ucp$ differ by a move, where
$a=\bc(\al,\be)$, $b=\bc(\be,\ga)$, and $b'=\bc(\be,\ga')$.  Since
$a$, $b$, and $b'$ pairwise share a common curve,
Proposition~\ref{sharethm} implies that $\ps(a)$, $\ps(b)$, and
$\ps(b')$ are pairwise sharing.

One can always find a separating curve $z$ which intersects $\ga'$ but
not any of $\al$, $\be$, or $\ga$ (by Lemma~\ref{4types}).  It follows
that $z$ intersects $b'$ but not $a$ or $b$, and then by
Lemma~\ref{disjlem}, we have that $\ps(z)$ intersects $\ps(b')$ but
not $\ps(a)$ or $\ps(b)$.

Suppose that $\ps(a)$, $\ps(b)$, and $\ps(b')$ do not all share the
same nonseparating curve.  Let $\pi\td\sigma\td\tau$ be a spine for
$\set{\ps(a),\ps(b)}$.  By the assumption, $\ps(b')$ does not share
$\sigma$ with $\ps(a)$ and $\ps(b)$, and hence it shares curves
$\omega$ and $\nu$ with $\ps(a)$ and $\ps(b)$, respectively.  Note
that $\sigma$, $\omega$, and $\nu$ must all be distinct, because
otherwise it follows that $\ps(b')$ is equal to either $\ps(a)$ or
$\ps(b)$.

We will now argue that there is no curve which intersects $\ps(b')$
and is disjoint from both $\ps(a)$ and $\ps(b)$.  This will contradict
our earlier statement about $\ps(z)$.  Indeed, any curve $c$ which
intersects $\ps(b')$ must also intersect at least one of $\omega$ or
$\nu$, say $\omega$.  Since $\omega$ lies on the genus 1 subsurface
bounded by $\ps(a)$, it follows that $c$ also intersects $\ps(a)$.
\epf

By Lemma~\ref{sharingpreserved}, well-definedness of $\pst$ is reduced
to the following proposition, which we will now deduce from general
work of Harer \cite{jlh}.

\begin{minithm}\label{connectivity}Any two spines $\al\td\be\td\ga$ and $\de\td\be\td\ep$ differ by a finite sequence of moves.\end{minithm}

\p{Harer's complex} Let $F$ be a surface with boundary, let $P$ be a
finite set of points in $\partial F$, and let $P_0$ be a subset of
$P$.  Harer defines an abstract simplicial complex $X=X(F,P,P_0)$
with:

{\bf Vertices}\qua Isotopy classes of arcs in $F$ connecting points of
$P_0$ to points of $P-P_0$.

{\bf Edges}\qua Two vertices are connected with an edge if the
corresponding arcs are disjoint (apart from endpoints) and their union
does not separate $F$.

In general, a $k$--simplex of $X$ is the span of $k+1$ pairwise
connected edges with the property that the union of the corresponding
arcs does not separate $F$.

We remark that Harer's original notation for this complex is $BX$.
Let $r'$ denote the number of boundary components of $F$ which contain
a point of $P$.  Harer proves:

\begin{thm}\label{harer}$X$ is spherical of dimension $2g-2+r'$.\end{thm}

By an argument of Hatcher \cite{ah}, we also have:

\begin{minithm}\label{chainconn}$X$ is chain-connected.  That is, any two maximal simplices are connected by a finite sequence of maximal simplices, where consecutive simplices in the sequence share a simplex of codimension 1.\end{minithm}

The key fact for this proposition is that the link of every simplex of
codimension greater than 1 is connected, as each such link is again
the complex $X$ for a surface with $2g-2+r' > 0$ (obtained by cutting
along the arcs represented by the simplex).  Thus, given any path
between maximal simplices, one can push the path off all simplices of
codimension greater than 1.  The new path gives the sequence of
maximal simplices in the statement of the proposition.

We now apply Harer's ideas to the current situation.  For $\be$ a
nonseparating curve, let $F$ be the surface obtained by cutting $S$
along $\be$.  Let $P=\set{p,q}$ be a pair of points, one on each
boundary component of $F$, and let $P_0=\set{p}$.  Vertices of the
complex $X=X(F,P,P_0)$ correspond to curves which intersect $\be$
once.  Further, every edge of $X$ is a spine of a sharing pair for
$\be$ (a pair of vertices fails to give a sharing pair only when the
corresponding arcs lie in a 3--holed sphere containing $\partial F$;
but in this case, the arcs necessarily separate $F$, so the vertices
do not form an edge).  Thus, a move between spines is achieved by
changing an edge to a new edge which lies in a triangle with the old
edge.

\bpf[Proof of Proposition~\ref{connectivity}]

Let $\uc$ and $\de\td\be\td\ep$ be two spines for the nonseparating
curve $\be$.  We think of these spines as edges in the complex
$X=X(F,P,P_0)$ defined above.  Let $M$ and $N$ be any maximal
simplices which contain these edges.  By Proposition~\ref{chainconn},
there is a sequence of maximal simplices of $X$:
\[ M = M_0, M_1, \dots, M_k = N \]
where $M_i$ and $M_{i+1}$ share a codimension 1 face.  Using this
sequence, we will construct the desired sequence of moves.  Let $e_0$
be the edge corresponding to $\uc$.  Now, assuming $e_i$ is an edge in
$M_i$ connecting $v$ to $w$, inductively define $e_{i+1}$ as follows:

If both $v$ and $w$ are vertices of $M_{i+1}$, then $e_{i+1}$ is
defined to be $e_i$.

In the case that $w$, say, is not a vertex of $M_{i+1}$ (note then
that $v$ must be in $M_{i+1}$), we define $e_{i+1}$ to be the span of
$v$ with any other vertex $w'$ of $M_{i+1} \cap M_i$.  Since $v$, $w$,
and $w'$ all lie in $M_i$, they form a triangle, and hence the edges
$vw$ and $vw'$ differ by a move.

Finally, $e_k$ differs from the edge corresponding to
$\de\td\be\td\ep$ by at most two moves, since they both lie in the
simplex $N$.
\epf

\subsection{Homeomorphism}\label{homeo}

At this point, we have shown that any superinjective map $\ps$ of
$\csep$ naturally gives rise to a map $\pst$ of $\ccs$.  In order to
show that $\pst$ is induced by a homeomorphism, we now check that
$\pst$ is a superinjective map of $\ccs$.  That is, for any curves $c$
and $d$, we show $\gin(c,d) \neq 0$ if and only if
$\gin(\pst(c),\pst(d)) \neq 0$:
\begin{enumerate}
\item If $c$ and $d$ are both separating, then apply the fact that
$\pst|_{\csep}=\ps$ is superinjective (Lemma~\ref{disjlem}).
\item If $c$ and $d$ are both nonseparating, then it suffices to note
that $c$ and $d$ are disjoint if and only if there are disjoint
sharing pairs representing those curves.
\item If $c$ is separating, and $d$ is nonseparating, then $c$ and $d$
are disjoint if and only if either $c$ is part of a sharing pair for
$d$ (if $d$ is on a genus 1 side of $c$) or $d$ has a sharing pair
whose curves are disjoint from $c$. \end{enumerate}

We appeal to the following theorem of Irmak \cite{ei}:

\begin{thm}\label{superinj}For $S$ a surface of genus $g \geq 3$, any superinjective map of $\ccs$ is induced by an element of $\mods$.\end{thm}

Thus, we have that $\pst$ is induced by $f \in \mods$, which simply
means that both give the same action on $\ccs$.

\section{Proofs of main theorems}\label{mainproof}

In Section~\ref{twistsec}, we showed that given any injection $\phi \co 
\fin \to \kg$ for $\fin$ a finite index subgroup of $\kg$, there is an
associated superinjective map $\ps$ of $\csep$.  Then, in
Section~\ref{nonsepsec}, we proved that any superinjective map $\ps$
of $\csep$ gives rise to a superinjective map $\pst$ of $\ccs$, which
is itself induced by a mapping class $f$.  Thus, we have:
\[ \phi \mapsto \ps \mapsto \pst \mapsto f \]

It is immediate that $\ps$ is induced by $f$, and hence we have proven
Theorem~\ref{superthm}.  As for the original map $\phi$, we now have
that any power of a twist $T_c^k \in \fin$ satisfies: \[ \phi(T_c^k) =
T_{f(c)}^{k'} \] for some $k' \in \z$, where $f$ is the element of
$\mods$ given by Theorem~\ref{superinj}.

We now show that $f$ induces $\phi$, by which we mean that for any $h
\in \fin$, its image is given by: \[ \phi(h) = fhf^{-1} \] From this,
Theorems~\ref{autsep} and~\ref{superthm}, as well as Main
Theorems~\ref{comm} and~\ref{intorelli}, will all follow.

Let $h \in \fin$.  Since two mapping classes which agree on all
separating curves are equal, it suffices to show that for any
separating curve $c$, we have $fh(c)=\phi(h)f(c)$.  Choose $k$ so that
$T_c^k,T_{h(c)}^k \in \fin$.  We have:
\begin{eqnarray*}T_{f \circ h(c)}^{k'} &=& \phi(T_{h(c)}^k) = \phi(hT_c^kh^{-1})  = \phi(h) \phi(T_c^k) \phi(h)^{-1} = \phi(h) T_{f(c)}^{k''} \phi(h)^{-1} \\ &=& T_{(\phi(h) \circ f)(c)}^{k''}  \end{eqnarray*}
But since Dehn twists can only be equal if they are twists about the
same curve, we have:
\[ (f \circ h) (c) = (\phi(h) \circ f)(c) \]
So it is true that $fh=\phi(h)f$, and hence $\phi(h) = fhf^{-1}$.

We have thus proven Main Theorem~\ref{intorelli}, which immediately
implies that $\phi$ is the restriction to $\fin$ of an element of
$\aut(\kg)$.  In particular: \[ \comm(\kg) \cong \aut(\kg) \]

Main Theorem~\ref{comm} now follows.  Indeed, the natural map
$\eta\co \mods \to \aut(\kg)$ is readily seen to be a right inverse of
the composition of maps taking $\phi$ to $f$, and thus $\mods \cong
\aut(\kg)$.  Theorem~\ref{autsep} follows similarly.

\section{Torelli group}\label{reproof}

We will now sketch how our theorems give alternate proofs of Farb and
Ivanov's results about Torelli groups.

Let $\psi\co  \fin \to \T$ be an injection, where $\fin$ is any finite
index subgroup of $\T$.  As in Farb and Ivanov's paper, $\psi$ induces
a superinjective map $\psi_\star$ of $\tg$ (this follows from
Proposition~\ref{trank}).  In particular, $\psi_\star$ induces a
superinjective map $\ps$ of $\csep$ (by the proof of
Lemma~\ref{genuslem}), and, following our paper, this gives rise to a
mapping class $f$ which induces $\ps$.  In summary: \[ \psi \mapsto
\psi_\star \mapsto \ps \mapsto f \] We now claim that $f$ also induces
$\psi_\star$.  It suffices to check $(f(a),f(b)) = \psi_\star((a,b))$
for any bounding pair $(a,b)$.  Given $(a,b)$ choose four separating
curves $w$, $x$, $y$, and $z$ as in Figure~\ref{bppic}.

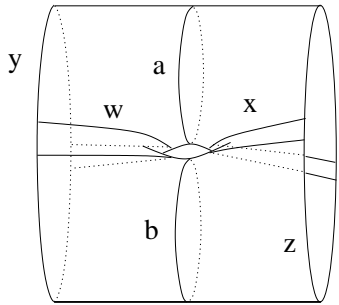
\pics{bp}{Separating curves defining a bounding pair}{4}

We have that $(f(a),f(b))$ is the unique bounding pair between $f(y)$
and $f(z)$ and disjoint from $f(w)$ and $f(x)$.

Now $\psi_\star$ has the same action as $f$ on $w$, $x$, $y$, and $z$
since all the curves are separating.  As in Section~\ref{basic},
$\psi_\star((a,b))$ must lie in between $f(y)$ and $f(z)$, and it must
be disjoint from $f(w)$ and $f(x)$.  This implies that
$\psi_\star((a,b))=(f(a),f(b))$.

By the same argument as in Section~\ref{mainproof}, $f$ induces
$\psi$.  Theorems~\ref{auttg} and~\ref{supertg} now follow, and from
these it is straightforward to prove Theorems~\ref{auttor}
and~\ref{cohopftor}.


\begin{thebibliography}

\bibitem{bm}
\textbf{Robert~W Bell}, \textbf{Dan Margalit}, \emph{Braid groups are almost
  co-{H}opfian}, \arxiv{math.GT/0403145}

\bibitem{blm}
\textbf{Joan~S Birman}, \textbf{Alex Lubotzky}, \textbf{John McCarthy},
  \emph{Abelian and solvable subgroups of the mapping class groups}, Duke Math.
  J. 50 (1983) 1107--1120 \MR{0726319}

\bibitem{bifa}
\textbf{Daniel Biss}, \textbf{Benson Farb}, \emph{${K}_g$ is not finitely
  generated},\nl \arxiv{math.GT/0405386}

\bibitem{bv}
\textbf{Martin~R Bridson}, \textbf{Karen Vogtmann}, \emph{Automorphisms of
  automorphism groups of free groups}, J. Algebra 229 (2000) 785--792
\MR{1769698}

\bibitem{fm} \textbf{Benson Farb}, \textbf{Lee Mosher}, \emph{The
geometry of surface-by-free groups}, Geom. Funct. Anal. 12 (2002)
915--963
\MR{1937831}

\bibitem{bf}
\textbf{Benson Farb}, \emph{{Automorphisms of the {T}orelli group}}, AMS
  sectional meeting, {A}nn {A}rbor, {M}ichigan, {M}arch 1, 2002

\bibitem{fh}
\textbf{Benson Farb}, \textbf{Michael Handel}, \emph{Commensurations of
  {O}ut(${F}_n$). {P}reprint, personal communication, {S}eptember 2004}

\bibitem{fi}
\textbf{Benson Farb}, \textbf{Nikolai~V Ivanov}, \emph{The {T}orelli geometry
  and its applications}, \arxiv{math.GT/0311123}

\bibitem{jlh}
\textbf{John~L Harer}, \emph{Stability of the homology of the mapping class
  groups of orientable surfaces}, Ann. of Math. 121 (1985) 215--249
\MR{0786348}

\bibitem{wjh}
\textbf{W\,J Harvey}, \emph{Boundary structure of the modular group}, from:
  ``Riemann surfaces and related topics (Stony Brook, 1978)'', Ann. of Math.
  Stud. 97, Princeton Univ. Press, Princeton, N.J. (1981)  245--251
\MR{0624817}


\bibitem{ah}
\textbf{Allen Hatcher}, \emph{On triangulations of surfaces}, Topology Appl. 40
  (1991) 189--194
\MR{1123262}

\bibitem{eiII} \textbf{Elmas Irmak}, \emph{Superinjective Simplicial
Maps of Complexes of Curves and Injective Homomorphisms of Subgroups
of Mapping Class Groups {II}}, \arxiv{math.GT/0311407}

\bibitem{ei}
\textbf{Elmas Irmak}, \emph{Superinjective simplicial maps of complexes of
  curves and injective homomorphisms of subgroups of mapping class groups},
  Topology 43 (2004) 513--541
\MR{2041629}

\bibitem{ninotes} \textbf{Nikolai~V Ivanov}, \emph{Automorphisms of
{T}eichm\"uller modular groups}, from: ``Topology and
geometry---Rohlin Seminar'', Lecture Notes in Math.  1346, Springer,
Berlin (1988) 199--270
\MR{0970079}

\bibitem{ni}
\textbf{Nikolai~V Ivanov}, \emph{Automorphism of complexes of curves and of
  {T}eichm\"uller spaces}, Internat. Math. Res. Notices  (1997) 651--666
\MR{1460387}

\bibitem{im}
\textbf{Nikolai~V Ivanov}, \textbf{John~D McCarthy}, \emph{On injective
  homomorphisms between {T}eichm\"uller modular groups. {I}}, Invent. Math. 135
  (1999) 425--486
\MR{1666775}

\bibitem{djI}
\textbf{Dennis Johnson}, \emph{The structure of the {T}orelli group. {I}. {A}
  finite set of generators for {${\mathcal I}$}}, Ann. of Math. 118 (1983)
  423--442
\MR{0727699}

\bibitem{djII}
\textbf{Dennis Johnson}, \emph{The structure of the {T}orelli group. {II}. {A}
  characterization of the group generated by twists on bounding curves},
  Topology 24 (1985) 113--126
\MR{0793178}

\bibitem{mk}
\textbf{Mustafa Korkmaz}, \emph{Automorphisms of complexes of curves on
  punctured spheres and on punctured tori}, Topology Appl. 95 (1999) 85--111
\MR{1696431}

\bibitem{fl}
\textbf{Feng Luo}, \emph{Automorphisms of the complex of curves}, Topology 39
  (2000) 283--298
\MR{1722024}

\bibitem{ms}
\textbf{Howard Masur}, \textbf{Saul Schleimer}, \emph{The Pants Complex Has
  Only One End}, \arxiv{math.GT/0312385}

\bibitem{mv}
\textbf{John~D McCarthy}, \textbf{William~R Vautaw}, \emph{Automorphisms of
  {T}orelli groups}, \arxiv{math.GT/0311250}

\bibitem{mm}
\textbf{Darryl McCullough}, \textbf{Andy Miller}, \emph{The genus {$2$}
  {T}orelli group is not finitely generated}, Topology Appl. 22 (1986) 43--49
\MR{0831180}

\bibitem{sm}
\textbf{Shigeyuki Morita}, \emph{Casson's invariant for homology {$3$}-spheres
  and characteristic classes of surface bundles. {I}}, Topology 28 (1989)
  305--323
\MR{1014464}

\bibitem{gp}
\textbf{Gopal Prasad}, \emph{Discrete subgroups isomorphic to lattices in
  semisimple {L}ie groups}, Amer. J. Math. 98 (1976) 241--261
\MR{0399351}

\bibitem{zs}
\textbf{Z Sela}, \emph{Structure and rigidity in ({G}romov) hyperbolic groups
  and discrete groups in rank {$1$} {L}ie groups. {II}}, Geom. Funct. Anal. 7
  (1997) 561--593
\MR{1466338}


\end{thebibliography}
\end{document}